\documentclass[12pt,a4paper]{article}
\usepackage{graphicx}
\begin{document}
 
 \thispagestyle{empty}
 
 \title{Platonic polyhedra tune the 3-sphere:\\ II. Harmonic analysis on  cubic spherical 3-manifolds.}
 \author{Peter Kramer, Institut fuer Theoretische Physik,\\ University  Tuebingen, Germany.}
 \maketitle

\section*{Abstract.}
From the homotopy groups of two distinct cubic spherical 3-manifolds we construct 
the isomorphic groups of deck transformations acting on the 3-sphere. These groups become the cyclic group of order 
eight and the quaternion group respectively.
By reduction of representations from the 
orthogonal group to the identity representation of these subgroups we provide two subgroup-periodic 
bases for the harmonic analysis on the 3-manifolds, which have applications to cosmic topology.

\section{Introduction.}
We view a spherical topological 3-manifold ${\cal M}$, see \cite{TH97}, as a prototile  on its  
cover $\tilde{{\cal M}}=S^3$. We study in \cite{KR08} the isometric actions of $O(4,R)$ on the 3-sphere $S^3$ and give its 
basis as  well-known homogeneous polynomials in \cite{KR05} eq.(37). From the homotopies of ${\cal M}$
we find as isomorphic images, see \cite{SE34}, its deck transformations on $S^3$. They form a subgroup 
$H$ of $O(4,R)$.  
By group/subgroup representation theory with intermediate Coxeter groups we construct on $S^3$ a $H$-periodic basis for the harmonic analysis on ${\cal M}$. 

Our approach yields in closed analytic form the onset, the selection rules, the 
multiplicity, projection operators, orthogonality rules and basis for each manifold. 
Among the Platonic polyhedra it  was applied  in \cite{KR05}, \cite{KR06}   to 
the Poincare dodecahedron  with $H$ the binary icosahedral group. 
An  algorithm due to Everitt in \cite{EV04} describes homotopies for   spherical 
3-manifolds from five Platonic polyhedra. Following it we found and applied in \cite{KR08}
for the tetrahedron  as $H$ the cyclic group $C_5$. 

One field of applications for  harmonic analysis  is cosmic topology, see \cite{LA95}, \cite{LU03},
\cite{AU05a}, \cite{AU05b}, \cite{AU08}: 
The topology of a 3-manifold ${\cal M}$ is favoured if 
data from the Cosmic Microwave Background can be  
expanded in its harmonic basis. 

Our work is well suited for comparing on $S^3$ 
a family of topologies and their  harmonic analysis. 
Here we turn to two distinct cubic spherical 3-manifolds,
with homotopies described in \cite{EV04}.
We employ an intermediate Coxeter group, construct the groups of deck transformations,  and derive 
and compare their harmonic analysis.

\section{The Coxeter group $G$ and the $8$-cell on $S^3$.}

The cartesian coordinates $x= (x_0,x_1,x_2,x_3)\in E^4$ for $S^3$ we combine   as in 
\cite{KR05}, \cite{KR08} in the matrix form
\begin{equation}
 \label{c3a}
u=
\left[
\begin{array}{ll}
 z_1&z_2\\
-\overline{z}_2&\overline{z}_1\\
\end{array}
\right],\:
z_1=x_0-ix_3,\; z_2=-x_2-ix_1,\: z_1\overline{z}_1+z_2\overline{z}_2=1.
\end{equation}
For the group action we start from the Coxeter group $G<O(4,R)$ \cite{HU90}, \cite{EV04} p. 254,
with the diagram
\begin{equation}
\label{c1}
G= \circ\stackrel{4}{-}\circ\stackrel{3}{-}\circ\stackrel{3}{-}\circ.
\end{equation}

\begin{center}
\includegraphics{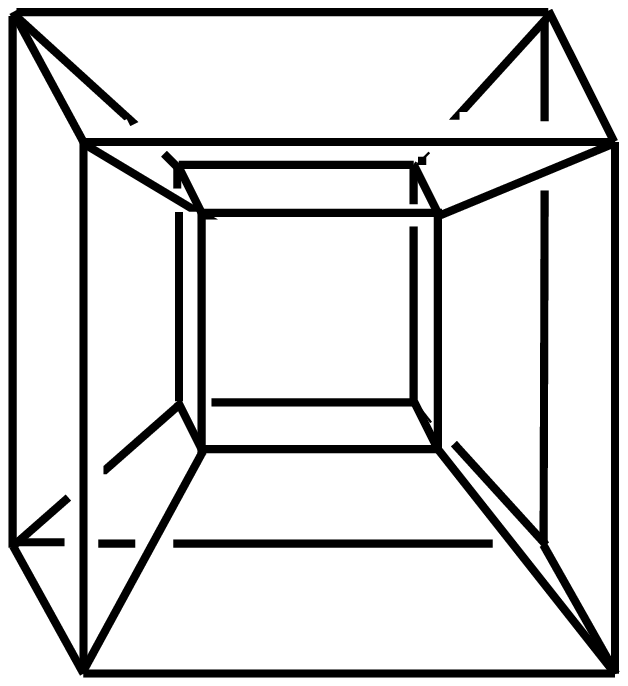}
\end{center}

Fig. 1. The 8-cell projected to the plane according to \cite{SO58} pp. 177-8. Any face bounds two, any edge  three, any vertex 
four out of the eight cubes.
\vspace{0.3cm}

For the Coxeter diagram  eq. \ref{c1} we  give for the $4$ Weyl reflections 
$W_s=W_{a_s}, s=1,2,3,4$
the Weyl vectors in {\bf Table 2.1} and compute for each Weyl vector 
$a_s=(a_{s0},a_{s1},a_{s2},a_{s3})$ 
the matrix 
\begin{equation}
\label{c2a}
v_s:=
\left[
\begin{array}{ll}
a_{s0}-ia_{s3}&-a_{s2}-ia_{s1}\\
a_{s2}-ia_{s1}&a_{s0}+ia_{s3}
\end{array} \right] \in SU(2,R).
\end{equation}

The matrices eq. \ref{c2a} will be used  to relate, see  \cite{KR08}, the 
Weyl reflections to  $SU^l(2,R)\times SU^r(2,R)$ acting by left and right multiplication 
on the coordinates eq. \ref{c3a}. 
We include in {\bf Table 2.1} the inversion ${\cal J}_4\in G$,  
and the additional Weyl operator $W_0$.

$G$ is isomorphic to the hyperoctahedral group, see \cite{CO65} p.90, which is a semidirect product of
its normal subgroup $(C_2)^4$ of four reflections in the four coordinate axes  with
the subgroup $S(4)$ of all $4!= 24$ permutations of these axes, 
\begin{equation}
 \label{c2}
G= (C_2)^4 \times_{s} S(4).
\end{equation}
The group is of order
$|G|=2^4 4!=384$.
Its elements we denote as products $g= \epsilon p, \epsilon \in (C_2)^4,\: p \in S(4)$
with the multiplication law 
\begin{eqnarray}
 \label{r0}
&&\epsilon= (\epsilon_0,\epsilon_1,\epsilon_2,\epsilon_3),\; \epsilon_j=\pm 1,
\\ \nonumber
&&\epsilon p \:\epsilon' p'= \epsilon'' p'',\;
\\ \nonumber
&&\epsilon''=\epsilon (p\epsilon'p^{-1}),\:p''=pp',
\\ \nonumber
&&(p\epsilon p^{-1})=
(\epsilon_{p^{-1}(0)},\epsilon_{p^{-1}(1)},\epsilon_{p^{-1}(2)},\epsilon_{p^{-1}(3)}).
\end{eqnarray}
The hyperoctahedral form of $G$ allows to study the irreducible representations.
We write the permutation $p\in S(4)$ in cycle form but use the numbers $\{0,1,2,3\}$ adapted to
the enumeration of the coordinates in eq. \ref{c3a}.
In {\bf Table 2.1} we give  the action of the Weyl reflections on the coordinates 
$x$ and the factors in $g=\epsilon p$.

\begin{center}
\includegraphics{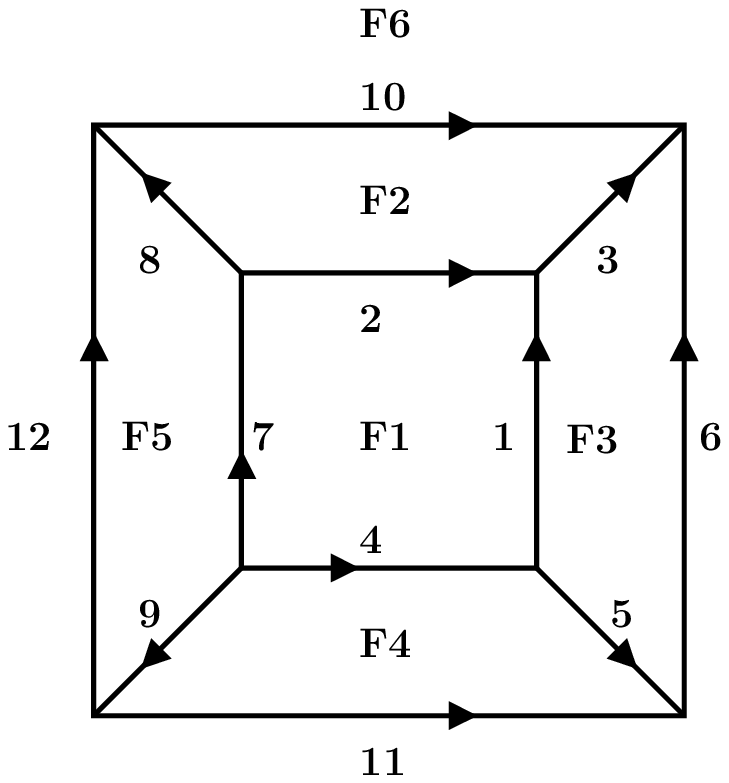} 
\end{center}

Fig. 2. Enumeration of faces $F1,\ldots, F6$ and edges $1,\ldots, 12$ for the
cubic prototile according to Everitt \cite{EV04} p. 260 Fig. 2.
\vspace{0.3cm}

The first three Weyl reflections from {\bf Table 2.1} generate, see \cite{HU90},  
the  cubic Coxeter subgroup 
\begin{equation}
 \label{c3}
O= \circ\stackrel{4}{-}\circ\stackrel{3}{-}\circ,
\end{equation}
isomorphic to the octahedral group $O \sim (C_2)^3 \times_{s} S(3)$
acting on $E^3\in E^4$. This group  from $48$ simplices  generates a spherical cube 
attached to a single vertex. We choose this cube  as the prototile on $S^3$.

We locate the center of the cubic prototile at $x=(1,0,0,0) \in S^3$.
When we include the action of $W_4$ along with the cubic Coxeter subgroup eq. \ref{c3},
this cube is mapped into $7$ companions which tile 
$S^3$ and generate the $8$-cell tiling described in \cite{SO58} pp. 177-8
and shown in a projection in Fig. 1.
The center positions of the $8$ spherical cubes are located 
at the $8$ points
\begin{equation}
\label{c6}
(\pm 1,0,0,0),(0,\pm 1,0,0),(0,0,\pm 1,0),(0,0,0,\pm 1). 
\end{equation}

Everitt \cite{EV04} has shown that $S^3$ admits two 
cubic spherical 3-manifolds with inequivalent first homotopy. 
He enumerates the faces and edges and  gives a graphical algorithm for their gluing, see 
Fig. 2. 
This gluing determines the generators for the first homotopy group.

\begin{center}
\includegraphics{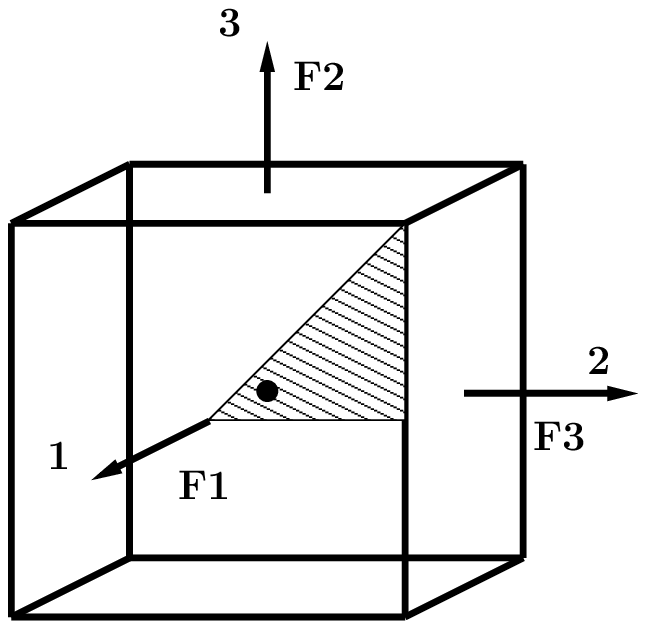} 
\end{center}

Fig. 3. Axes ${\bf (x_1,x_2,x_3)}\in E^3$ marked by arrows for the  cubic prototile with face enumeration $(F1, F3, F2)$
according to Fig. 2. The three axes are cyclically rotated into one another 
by the rotation $(W_2W_3),\; (W_2W_3)^3=e$.
\vspace{0.3cm}

In Fig. 3, the three edges of the shaded triangle mark the intersections of  the Weyl reflection (hyper-)planes 
for $W_1,W_2,W_3$ with the face $F1$ of the cube. Face $F1$ itself is part of the Weyl reflection 
(hyper-)plane for $W_4$. In the homotopy group, there appears a gluing of opposite faces.
Use of the inversion ${\cal J}_3$ in the center (black circle), followed by the Weyl reflection $W_4$, converts
this gluing  of opposite faces into a standard deck operation 
$st(1\Leftarrow 6)$ eq. \ref{c12}. This operation maps  
the ${\cal J}_3$-inverted and $W_4$-reflected initial cube into a new position with its face $F6$ 
glued to face $F1$ of the
initial cube.
\vspace{0.3cm}

{\bf Table 2.1}
The Weyl vectors $a_s,\: s=1,.., 4$ and $a_0$ for the Coxeter group $G$ eq. \ref{c1},  the
$2 \times 2$ unitary matrices $v_s$ eq. \ref{c2a}, the action on $x$, and the product form 
$g=\epsilon p$ with $p$ in cycle form.
\begin{equation}
\label{c7}
\begin{array}{l|l|l|l|l|l}
s &{\rm Weyl}\:{\rm vector}\:a_s &{\rm matrix}\: v_s&g\,x&\epsilon&p\\
\cline{1-6}
&&\\ 
1 & (0,0,0,1)
&\left[
\begin{array}{ll}
-i&0\\
0&i\\
\end{array}
\right]&
(x_0,x_1,x_2,-x_3)&(+++-)&e
\\
2 & (0,0,-\sqrt{\frac{1}{2}},\sqrt{\frac{1}{2}})
&\sqrt{\frac{1}{2}}
\left[\begin{array}{ll}
-i&1\\
-1&i\\
\end{array}
\right]&(x_0,x_1,x_3,x_2)&(++++)&(23)
\\
3 & (0,\sqrt{\frac{1}{2}},-\sqrt{\frac{1}{2}},0)
&\sqrt{\frac{1}{2}}
\left[
\begin{array}{ll}
0&1-i\\
-1-i&0\\
\end{array}
\right]&(x_0,x_2,x_1,x_3)&(++++)&(12)
\\
4 & (-\sqrt{\frac{1}{2}},\sqrt{\frac{1}{2}},0,0)
&\sqrt{\frac{1}{2}}
\left[
\begin{array}{ll}
-1&-i\\
-i&-1\\
\end{array}
\right]&(x_1,x_0,x_2,x_3)&(++++)&(01)
\\
0 & (1,0,0,0)
&\left[
\begin{array}{ll}
1&0\\
0&1\\
\end{array}
\right]&(-x_0,x_1,x_2,x_3)&(-+++)&e
\\
&{\cal J}_4&&(-x_0,-x_1,-x_2,-x_3)&(----)&e
\end{array}
\end{equation}

\section{The  group isomorphism  ${\rm deck}(C2) \sim \pi_1(C2) \sim C_8$ for the cubic spherical 3-manifold $(C2)$.}

We start from the algorithm of \cite{EV04} p. 259 Table 3 on  the first homotopy group $\pi_1(C2)$ for the cube and construct the explicit
isomorphism to the group of deck transformations ${\rm deck}(C2)$.
We denote the faces of the cube from Fig. 2 as $F1,\ldots, F6$. To the first three faces of $C2$,
the glue partners are
\begin{equation}
\label{c9}
F3\cup F1,\; F4\cup F2,\; F5\cup F3.
\end{equation} 
We represent the four edges for a square face by the corresponding numbers from Fig. 2. 
Evaluating the glue  algorithm  for $C2$  in \cite{EV04},  the first edge glue generator  can be depicted   from right to left as the map
\begin{equation}
\label{c10}
g_1=g_1(1\Leftarrow 3):
\\ \nonumber
\left[
\begin{array}{rrr}
&\overline{1}&\\
\overline{4}&&2\\
&7&\\
\end{array}
\right]
  \Leftarrow 
\left[
\begin{array}{rrr}
&\overline{3}&\\
\overline{1}&&6\\
&5&\\
\end{array}
\right].
\end{equation}
We use left arrows $\Leftarrow$ in line with the usual sequence of operator products.
Now we wish to map this glue generator $g_1(1\Leftarrow 3)$ isomorphically into a generator  
from the group  ${\rm deck}(C2)$ of deck transformation acting on the 3-sphere $S^3$.
We follow a similar method as employed  in \cite{KR05} pp. 5322-4 for the Poincare dodecahedral 
3-manifold.
We refer to Fig. 3   for the position of the cube and for the three orthogonal directions
$1,2,3$. We factorize $g_1(1\Leftarrow 3)$ into three actions: A positive rotation  
$R_3(\pi/2)$ around the 3-axis, followed by a standard glue operator $st(1\Leftarrow 6)$
eq. \ref{c12},  followed 
by a positive rotation $R_1(\pi/2)$ around the 1-axis. In eq. \ref{c11} 
we depict from right to left the images of the initial face $F3$. In the second line we 
write the corresponding triple product of operators.
\begin{eqnarray}
\label{c11}
&& g_1=g_1(1\Leftarrow 3):
\\ \nonumber
&&
\begin{array}{ccccccc}
\left[
\begin{array}{rrr}
&\overline{1}&\\
\overline{4}&&2\\
&7&\\
\end{array}
\right]
& \Leftarrow &
\left[
\begin{array}{rrr}
&\overline{4}&\\
7&&\overline{1}\\
&2&\\
\end{array}
\right]
& \Leftarrow &
\left[
\begin{array}{rrr}
&10&\\
\overline{6}&&12\\
&\overline{11}&\\
\end{array}
\right]
& \Leftarrow &
\left[
\begin{array}{rrr}
&\overline{3}&\\
\overline{1}&&6\\
&5&\\
\end{array}
\right]
\\
& R_1(\frac{\pi}{2}) 
&\times 
&{\rm st}(1\Leftarrow 6)
& \times 
&R_3(\frac{\pi}{2})&\\
\end{array}
\end{eqnarray}
The orientation of edges is counter-clockwise for the two right-hand diagrams but
clockwise for the left-hand ones.
The third operator maps $F6 \Leftarrow F3$, the first one rotates $F1$ while
preserving the center $x=(1,0,0,0)$ of the initial cube. The second operator is crucial for the isomorphism $\pi_1(C2)\rightarrow deck(C2)$.
It is constructed as the product 
\begin{equation}
\label{c12}
 {\rm st}(1\Leftarrow 6)=W_4{\cal J}_3,\; {\cal J}_3=W_0 {\cal J}_4,
\end{equation}
of the inversion  operator ${\cal J}_3 \in H$ w.r.t. $(x_1,x_2,x_3)$, which can be factorized 
into the Weyl operator $W_0$ and the full coordinate inversion ${\cal J}_4\in G$,
followed by the Weyl reflection $W_4$.
The operator ${\cal J}_3$ in eq. \ref{c11} inverts the cube in its center and so maps any face $F_i$ of the cube into 
its opposite face, in the present case $F1 \Leftarrow F6$. The final Weyl reflection 
$W_4$ following ${\cal J}_3$ reflects the inverted cube in its face $F1$ and produces 
an image  of the initial  cube with the new center $(0,1,0,0)$, glued with its face $F3$ to the original face $F1$. The product eq. \ref{c12} contains two Weyl reflections and therefore preserves
the orientation.

For later use we list three  standard glue operators of the cube in {\bf Table 3.1}.
We construct them from eq. \ref{c12} by the conjugations
\begin{equation}
 \label{c13}
{\rm st}(2\Leftarrow 4)=(W_3W_2){\rm st}(1\Leftarrow 6)(W_2W_3),\;
{\rm st}(3\Leftarrow 5)=(W_2W_3){\rm st}(1\Leftarrow 6)(W_3W_2).
\end{equation} 
listed in {\bf Table 3.2}.
{\bf Table 3.1}: Standard glue operators in the Coxeter group $G$.
\begin{equation}
 \label{c26}
 \begin{array}{l|lll}
g         &g\: x&\epsilon&p\\ \hline
st(1\Leftarrow 6) &(x_1,-x_0,-x_2,-x_3)&(+---)&(01)\\
st(2\Leftarrow 4)&(x_2,-x_1,-x_0,-x_3)&(+---)&(02)\\
st(3\Leftarrow 5)&(x_3,-x_1,-x_2,-x_0)&(+---)&(03)\\   
\end{array}
\end{equation}

The operations eq. \ref{c13} are in $G$ but not necessarily in ${\rm deck}(C2)$.
The rotations and the product $W_4W_0$ appearing in eq. \ref{c12} are even in $G$, and 
so all operations preserve  orientation.

Next we express the generator $g_1$ as a product of pairs of Weyl operators and find
\begin{eqnarray}
 \label{c14} 
&& R_1(\pi/2)=(W_2W_1),\: st(1\Leftarrow 6)= (W_4W_0){\cal J}_4,
\\ \nonumber 
&&R_3(\pi/2)=(W_2W_3)(W_2W_1)(W_3W_2),
\\ \nonumber
&&g_1=(W_2W_1)(W_4W_0){\cal J}_4(W_2W_3)(W_2W_1)(W_3W_2).
\end{eqnarray}
We use \cite{KR08} eqs. (60) for products of Weyl operators in eq. \ref{c14}  and the matrices $v_s$ from {\bf Table 2.1}
to rewrite the operator $T_{g_1}$ in the form
\begin{eqnarray}
 \label{c15}
&&T_{g_1}= T_{(w_{l1},w_{r1})},
\\ \nonumber
&&w_{l1}=(v_2v_1^{-1})(v_4v_0^{-1})(v_2v_3^{-1})
(v_2v_1^{-1})(v_3v_2^{-1})
=\sqrt{\frac{1}{2}}
\left[
\begin{array}{ll}
-\overline{a}&0\\
0&-a\\
\end{array}
\right],
\\ \nonumber
&&w_{r1}= (-1)(v_2^{-1}v_1)(v_4^{-1}v_0)(v_2^{-1}v_3)
(v_2^{-1}v_1)(v_3^{-1}v_2)=
\sqrt{\frac{1}{2}} 
\left[
\begin{array}{ll}
0&-a^3\\
-a& 0\\
\end{array}
\right],\: a=\exp(\pi i/4).
\end{eqnarray}
By $(w_{l},w_{r})$ we denote the elements of the $SU^l(2,R)\times SU^r(2,R)$ 
action $u \rightarrow w_{l}^{-1}uw_{r}$ on $S^3$ in the coordinates eq. \ref{c3a}. In this scheme, 
${\cal J}_4$ in eq. \ref{c12} yields  the operator
\begin{equation}
 \label{c16}
T_{{\cal J}_4}=T_{(e,-e)}.
\end{equation}
which commutes with all rotation operators. In eq. \ref{c15} we have absorbed this 
operator into $w_r$.
From eq. \ref{c15} we easily see that $w_{l1}^8=w_{r1}^8=e$ so that $g_1,\: g_1^8=e$ generates the cyclic group 
$C_8$ of order $8$.

We can construct from the graphs given in \cite{EV04} the other glue generators of the homotopy group of $C2$ and map them
isomorphically into generators of the group of deck transformations. It turns out that 
all of these glue generators   become powers of the first glue generator  $g_1 \in C_8$ eq.  and its image eq. \ref{c15}  in $deck(C2)$.
We list the actions of the  eight powers $g_1^t$ in {\bf Table 3.2}.

{\bf 1 Theorem}: The homotopy group and the group of deck transformations for the spherical cubic 3-manifold 
$C2$ of Everitt \cite{EV04} p. 259 Table 3 are isomorphic to the cyclic group 
\begin{equation}
 \label{c19}
{\rm deck}(C2)= C_8= \langle g_1^t,\: t=1, \ldots, 8,\: g_1^8=e\rangle
\end{equation}
with actions on $S^3$  given in {\bf Table 3.2}. The deck transformations   generate fix-point free  the 8-cell on $S^3$.
The group/subgroup scheme for the cubic 3-manifold $C2$ is
\begin{equation}
 \label{c20}
O(4,R)>G>C_8.
\end{equation}

\subsection{The reduction $O(4,R)>C_8$ and harmonic analysis on (C2).}
The irreducible representations of $O(4,R)$, their polynomial basis 
in terms of spherical harmonics, and the matrix elements of Weyl operators  
we adopt from \cite{KR08} section 4.2.
Here we consider directly the reduction of representations for the group/subgroup pair
\begin{equation}
 \label{c21}
O(4,R)>C_8.
\end{equation}
The operators for the generator of $g_1 \in C_8$ and its powers from \cite{KR08}
eq. (60)  have the form
\begin{equation}
\label{c22}
T_{g_1}= T_{(w_{l1},w_{r1})},\: T_{g_1^t}=T_{(w_{l1}^t,w_{r1}^t)},  
\end{equation}
given in {\bf Table 3.2} in terms of the $SU^l(2,R)\times SU^r(2,R)$ action.
\vspace{0.2cm}

{\bf Table 3.2}: The elements of the cyclic group $C_8$ of deck transformations  of
the manifold $C2$ and their actions on $S^3$.

\begin{equation}
\label{c28}
\begin{array}{l|lllll}
t& (g_1)^t\;x &w_l^t&w_r^t&\epsilon&p \\
\hline
1& (x_1,-x_3,x_0,x_2)&
\left[\begin{array}{ll}
-\overline{a}&0\\
0&-a\\
\end{array}
\right]
&
\left[\begin{array}{ll}
0&-a^3\\
-a&0\\
\end{array}
\right]
&(+-++)&(0132)\\
2& (-x_3,-x_2,x_1,x_0)&
\left[\begin{array}{ll}
\overline{a}^2&0\\
0&a^2\\
\end{array}
\right]
&
\left[\begin{array}{ll}
-1&0\\
0&-1\\
\end{array}
\right]
&(--++)&(03)(12)\\
3& (-x_2,-x_0,-x_3,x_1)&
\left[\begin{array}{ll}
-\overline{a}^3&0\\
0&-a^3\\
\end{array}
\right]
&
\left[\begin{array}{ll}
0&a^3\\
a&0\\
\end{array}
\right]
&(---+)&(0231)\\
4& (-x_0,-x_1,-x_2,-x_3)&
\left[\begin{array}{ll}
-1&0\\
0&-1\\
\end{array}
\right]
&
\left[\begin{array}{ll}
1&0\\
0&1\\
\end{array}
\right]
&(----)&e\\
5& (-x_1,x_3,-x_0,-x_2)&
\left[\begin{array}{ll}
\overline{a}&0\\
0&a\\
\end{array}
\right]
&
\left[\begin{array}{ll}
0&-a^3\\
-a&0\\
\end{array}
\right]
&(-+--)&(0132)\\
6& (x_3,x_2,-x_1,-x_0)&
\left[\begin{array}{ll}
-\overline{a}^2&0\\
0&-a^2\\
\end{array}
\right]
&
\left[\begin{array}{ll}
-1&0\\
0&-1\\
\end{array}
\right]
&(++--)&(03)(12)\\
7& (x_2,x_0,x_3,-x_1)&
\left[\begin{array}{ll}
\overline{a}^3&0\\
0&a^3\\
\end{array}
\right]
&
\left[\begin{array}{ll}
0&a^3\\
a&0\\
\end{array}
\right]
&(+++-)&(0231)\\
8& (x_0,x_1,x_2,x_3)&
\left[\begin{array}{ll}
1&0\\
0&1\\
\end{array}
\right]
&
\left[\begin{array}{ll}
1&0\\
0&1\\
\end{array}
\right]
&(++++)&e\\
\end{array}
\end{equation}

From eq. \ref{c22} and from \cite{KR08} eq. (45) we find for the characters of the powers $g_1^t$ 
in the irreducible representation $D^{(j,j)}$ in terms of characters of $SU
(2,R)$ the result 
\begin{equation}
 \label{c23}
\chi^{(j,j)}(g_1^t)= \chi^j(w_{l1}^t)\chi^j(w_{r1}^t),\; t=1,\ldots, 8.
\end{equation}
These expressions are easily evaluated with the help of \cite{KR08}, Appendix.
The multiplicity of the identity representation $D^0$ of $C_8$ is then given 
by 
\begin{equation}
\label{c24}
m(C_8 (j,j),0)= \frac{1}{8}\sum_{t=1}^8\;  \chi^j(w_{l1}^t)\chi^j(w_{r1}^t).
\end{equation}

In {\bf Table 3.3} we give the values of the characters of $SU^l(2,C),SU^^r(2,C)$ and the
onset of multiplicities of $C_8$-periodic states eq. \ref{c24} as functions of  the degree $2j,\: 0\leq j\leq 8$. 
The characters $\chi^j$ in this table can be divided into a set with the property 
$\chi^{j+4}(w_{l,r}^t)=\chi^j(w_{l,t}^t)$ and a set which increases with $j$. 
Treating these sets in eq. \ref{c24} separately allows to derive the 
following recursion relation for the multiplicities:
\begin{equation}
 \label{c25} 
(2j)= {\rm even}: m(C_8 (j+4,j+4),0) =m(C_8 (j,j),0)+8j+20+2(-1)^j. 
\end{equation}
Multiplicities for  odd degree $(2j)$ vanish: From {\bf Table 3.2}
one finds $g_1^4={\cal J}_4,\: {\cal J}_4 \in C_8$. The representation of any  group element applied to a $C_8$-periodic polynomial must give unity. For the inversion ${\cal J}_4$ in {\bf Table 2.1}  this allows only for even  
degree.
\vspace{0.2cm}

{\bf Table 3.3}: The characters $\chi^j(w_l^t), \chi^j(w_r^t)\: t=1,\ldots,8$ in $SU(2,R)$  for $C_8$ in the irreducible representations $D^j$ of $SU(2,R)$.
The characters  $\chi^{(j,j)}(g_1^t)$ are given by eq.\ref{c23} and the multiplicities 
$m(C_8(j,j)0)$ by eq. \ref{c24}.

\begin{eqnarray}
 \label{c28a}
&&\begin{array}{lll|rrrrrrrrr}
             &\phi/2&j:&           0&1&2 &3 &4 &5&6&7&8\\ \hline
\chi^j(w_l)  &3\pi/4&&1&1&-1&-1&1 &1 &-1&-1&1\\
\chi^j(w_l^2)&\pi/2&&1&-1&1&-1&1&-1&1&-1&1\\
\chi^j(w_l^3)&\pi/4&&1&1&-1&-1&1&1&-1&-1&1\\
\chi^j(w_l^4)&\pi&&1&3&5&7&9&11&13&15&17\\
\chi^j(w_l^5)&\pi/4&&1&1&-1&-1&1&1&-1&-1&1\\
\chi^j(w_l^6)&\pi/2&&1&-1&1&-1&1&-1&1&-1&1\\
\chi^j(w_l^7)&\pi/4&&1&1&-1&-1&1&1&-1&-1&1\\
\chi^j(w_l^8)&0&&1&3&5&7&9&11&13&15&17\\ \hline
\chi^j(w_r)  &\pi/2&&1&-1&1&-1&1 &-1 &1&-1&1\\
\chi^j(w_r^2)&\pi&&1&3&5&7&9&11&13&15&17\\
\chi^j(w_r^3)&\pi/2&&1&-1&1&-1&1&-1&1&-1&1\\
\chi^j(w_r^4)&0&&1&3&5&7&9&11&13&15&17\\
\chi^j(w_r^5)&\pi/2&&1&-1&1&-1&1&-1&1&-1&1\\
\chi^j(w_r^6)&\pi&&1&3&5&7&9&11&13&15&17\\
\chi^j(w_r^7)&\pi/2&&1&-1&1&-1&1&-1&1&-1&1\\
\chi^j(w_r^8)&0&&1&3&5&7&9&11&13&15&17\\ \hline
m(C_8 (j,j),0) &&&1&1&7&11&23&27&45&53&77\\
\end{array}
\\ \nonumber
\end{eqnarray}

\subsection{An orthogonal $C_8$-periodic basis for the harmonic analysis on $C2$.}

The spherical harmonics on $S^3$ are given from \cite{KR08} eq. (85) as the polynomials $D^j_{m_1,m_2}(z_1,z_2,\overline{z}_1,\overline{z}_2)$. 
Upon inserting the Euler angles $(\alpha,\beta,\gamma)$ by
\begin{equation}
 \label{s4}
z_1=\exp (i(\alpha+\gamma)/2)\cos(\beta/2)), z_2=\exp (i(\alpha-\gamma)/2)\sin(\beta/2)),
\end{equation}
into these spherical harmonics, we obtain the Wigner representation matrices of $u \in SU(2,C)$.
Here we need them only for integer $j$.
With the Euler angles we can then write the orthogonality relations from \cite{ED57} eq.(4.6.1) as 
\begin{equation}
 \label{s5}
\frac{1}{8\pi^2}\int_0^{2\pi}\int_0^{\pi}\int_0^{2\pi}
\overline{D}^{j'}_{m_1',m_2'}(\alpha \beta \gamma)
D^{j}_{m_1,m_2}(\alpha \beta \gamma)d\alpha \sin \beta d\beta d\gamma
=\delta_{j',j}\;\delta_{m_1',m_1}\;\delta_{m_2',m_2}\frac{1}{2j+1}.
\end{equation}
With the representations of $SU^l(2,C) \times SU^r(2,C)$ given in \cite{KR08} eq. (85)
we construct from {\bf Table 3.2} 
the representations of $C_8$ of the form
\begin{equation}
 \label{cc1}
(T_{(w_l^i,w_r^i)}D)^j_{m_1,m_2}(u)=D^j_{m_1,m_2}(w_l^{8-i}uw_r^i)= 
\sum_{m_1',m_2'} D^j_{m_1',m_2'}(u) D^j_{m_1,m_1'}(w_l^{8-i})D^j_{m_2',m_2}(w_r^i),
\end{equation}
whose coefficients are given in {\bf Table 3.4}.
\vspace{0.2cm}

{\bf Table 3.4}: Matrix elements in the representation $D^{(j,j)}$ eq. \ref{cc1} of the
elements of $C_8$ for $i=1,\ldots, 8$ with  $l=1,\ldots, 4$.

\begin{equation}
\begin{array}{l|ll} 
i&D^j_{m_1,m_1'}(w_l^{8-i})&D^j_{m_2',m_2}(w_r^i)\\
\hline
1& i^{m_1}\delta_{m_1,m_1'}&(-1)^{2j}(-1)^{j+m_2}i^{m_2}\delta_{m_2',-m_2}\\
2& (-1)^{2j}(-1)^{m_1}\delta_{m_1,m_1'}&(-1)^{2j}\delta_{m_2',m_2}\\
3& (-1)^{m_1} i^{m_1}\delta_{m_1',m_1}&(-1)^{j+m_2}i^{m_2}\delta_{m_2',-m_2}\\
4& (-1)^{2j}\delta_{m_1,m_1'}&\delta_{m_2',m_2}\\
4+l& (-1)^{2j}D^j_{m_1,m_1'}(w_l^{8-l})&D^j_{m_2',m_2}(w_r^l)\\
\end{array}
\end{equation}
\vspace{0.2cm}

The projection operator on the identity representation $D^0$ of $C_8$ is
\begin{equation}
\label{cc2}
P^{0}=\frac{1}{8}\sum_i^8 T_{(w_l^i,w_r^i)}.  
 \end{equation}
Its matrix elements from {Table 3.4} for fixed $j$ can be rewritten as
\begin{eqnarray}
 \label{cc3}
&& (P^{0,j})_{m_1,m_2,m_1',m_2'}
\\ \nonumber 
&&=\frac{1}{2}\left[1+(-1)^{2j}\right]\frac{1}{2}\left[1+(-1)^{m_1}\right]
\\ \nonumber 
&&\times \frac{1}{2}
\left[\delta_{m_1,m_1'}\delta_{m_2',m_2}
+i^{m_1}(-1)^{j+m_2}i^{m_2}\delta_{m_1,m_1'}\delta_{m_2',-m_2}\right]
\end{eqnarray}
The two prefactors imply integer $j$ and even $m_1$. For this 
case we get by application of the projector 
\newpage 

{\bf Table 3.5}: The $C_8$-periodic basis $\{\phi^j_{m_1,m_2}\}$ on $S^3$ for the harmonic analysis on
the cubic spherical 3-manifold $C2$ in terms of spherical harmonics $D^j(u)$ on $S^3$.
\begin{eqnarray}
 \label{cc4}
&&j={\rm integer},\, m_1={\rm even},\,-j\leq m_1\leq j,\, i^{m_1}(-1)^j=1, \, m_2=0:
\\ \nonumber 
&&\phi^j_{m_1,0}=\frac{\sqrt{2j+1}}{\sqrt{8}\pi}D^j_{m_1,0}(u),
\\ \nonumber
&&j={\rm integer},\, m_1={\rm even},\,-j\leq m_1\leq j,\, 0<m_2 \leq j:
\\ \nonumber
&&\phi^j_{m_1,m_2}=\frac{\sqrt{2j+1}}{4\pi}
\left[D^j_{m_1,m_2}(u)+i^{m_1}(-1)^{(j+m_2)}i^{m_2}D^j_{m_1,-m_2}(u)\right] 
\end{eqnarray}
{\bf 3 Theorem}: An orthonormal basis for the harmonic analysis on the cubic spherical manifold 
$C2$ is spanned by the $C_8$-periodic polynomials of degree $2j, j=0,1,2,...$ in {\bf Table 3.5}.

\section{The  group isomorphism  $\pi_1(C3) \sim {\rm deck}(C3)\sim Q$ for  the cubic spherical 3-manifold $(C3)$.}
The second  possible homotopy group of the spherical 3-cube is given by Everitt \cite{EV04} 
in Table 3 p. 259 by a second graphical algorithm. We denote this cubic  3-manifold 
as $C3$.  The order of the homotopy group and the group of deck transformations is  again $8$. 
For the homotopy of $C3$, 
opposite faces of the cube are glued,
\begin{equation}
\label{c34} 
F6\cup F1,\, F5\cup F3,\: F4\cup F2.
\end{equation}
From the edge gluing we find that each  glue is followed by  a left-handed rotation by $\pi/2$.  We construct three glue generators $q_1, q_2, q_3$ from the prescription of 
\cite{EV04} p. 259 Table 3, depict the edge gluings from right to left  as before
and factorize them with the standard glue operators eq. \ref{c13}:
\begin{eqnarray}
\label{c35} 
&&q_1=q_1(1\Leftarrow 6):
\\ \nonumber
&&
\begin{array}{ccccc}
\left[
\begin{array}{rrr}
&7&\\
2&&\overline{4}\\
&\overline{1}&\\
\end{array}
\right]
&  \Leftarrow 
&
\left[
\begin{array}{rrr}
&10&\\
\overline{6}&&12\\
&\overline{11}&\\
\end{array}
\right],
&&\\
\left[
\begin{array}{rrr}
&7&\\
2&&\overline{4}\\
&\overline{1}&\\
\end{array}
\right]
&  \Leftarrow 
&
\left[
\begin{array}{rrr}
&\overline{4}&\\
7&&\overline{1}\\
&2&\\
\end{array}
\right]
&  \Leftarrow 
&
\left[
\begin{array}{rrr}
&10&\\
\overline{6}&&12\\
&\overline{11}&\\
\end{array}
\right]
\\
&R_1(-\pi/2)&\times &st(1\Leftarrow 6)&\\
\end{array}
\\ \nonumber
&&q_2=q_2(3\Leftarrow 5):
\\ \nonumber
&&
\begin{array}{ccccc}
\left[
\begin{array}{rrr}
&1&\\
3&&\overline{5}\\
&\overline{6}&\\
\end{array}
\right]
&  \Leftarrow 
&
\left[
\begin{array}{rrr}
&8&\\
\overline{12}&&7\\
&\overline{9}&\\
\end{array}
\right],
&&\\
\left[
\begin{array}{rrr}
&1&\\
3&&\overline{5}\\
&\overline{6}&\\
\end{array}
\right]
&  \Leftarrow 
&
\left[
\begin{array}{rrr}
&\overline{5}&\\
1&&\overline{6}\\
&3&\\
\end{array}
\right]
&  \Leftarrow 
&
\left[
\begin{array}{rrr}
&8&\\
\overline{12}&&7\\
&\overline{9}&\\
\end{array}
\right]
\\
&R_2(-\pi/2)&\times&st(3\Leftarrow 5)&\\
\end{array}
\\ \nonumber
&&q_3=q_3(2\Leftarrow 4):
\\ \nonumber
&&
\begin{array}{ccccc}
\left[
\begin{array}{rrr}
&\overline{3}&\\
\overline{2}&&10\\
&8&\\
\end{array}
\right]
&  \Leftarrow 
&
\left[
\begin{array}{rrr}
&4&\\
9&&\overline{5}\\
&11&\\
\end{array}
\right],
&&\\
\left[
\begin{array}{rrr}
&\overline{3}&\\
\overline{2}&&10\\
&8&\\
\end{array}
\right]
&  \Leftarrow 
&
\left[
\begin{array}{rrr}
&10&\\
\overline{3}&&8\\
&\overline{2}&\\
\end{array}
\right]
&  \Leftarrow 
&
\left[
\begin{array}{rrr}
&\overline{4}&\\
9&&\overline{5}\\
&11&\\
\end{array}
\right]
\\
&R_3(-\pi/2)&\times&st(2\Leftarrow 4)&\\
\end{array}
\end{eqnarray}
The three generators in terms of Weyl reflections become
\begin{eqnarray}
\label{c35a}
&&q_1=(W_1W_2)(W_4W_0){\cal J}_4,
\\ \nonumber 
&&q_2=(W_3W_2)q_1(W_2W_3),
\\ \nonumber 
&&q_3=(W_2W_3)q_1(W_3W_2).
\end{eqnarray}
Analysis of the products of Weyl operators in eq. \ref{c35a} in terms of actions
of $SU^l(2,R)\times SU^r(2,R)$ given in eqs. \ref{c42}, \ref{c43} show that the generators 
$q_1,q_2,q_3$ act exclusively from the left on $u \in S^3$ eq. \ref{c3a} and in $SU^l(2,R)$ 
take a form  equivalent according to 
\begin{equation}
 \label{c35b}
(q_1,q_2,q_3)\sim (-{\bf k},-{\bf j},-{\bf i})
\end{equation}
with  the standard $2 \times 2$ quaternion representation. 
 
As mentioned before,
any apparent reversion of the edge orientation is compensated by a Weyl reflection.
Expressed in $G$ the generators from eq. \ref{c35a} are  given in  {\bf Table 4.1}.
Their multiplication yields the relations
\begin{equation}
\label{c36}
q_1^2=q_2^2=q_3^2=q_1q_2q_3= {\cal J}_4. 
\end{equation}
These are exactly the relations characterizing the quaternion group $Q$, see \cite{CO65} p.8,
with ${\cal J}_4$ playing the role of $(-1)$.
The eight elements are
\begin{equation}
\label{c37}
{\rm deck}(C3)= Q = \langle e,  q_j,{\cal J}_4, {\cal J}_4q_j,\: j=1,2,3\rangle. 
\end{equation}

We therefore have shown\\ 
{\bf 2 Theorem}: The homotopy group and the group of deck transformations for the spherical cubic 3-manifold 
$C3$ of Everitt \cite{EV04} p. 259 Table 3 are isomorphic to the quaternion group $Q$. The elements 
of the group
${\rm deck}(C3)$ act on  $S^3$ from the left as given in eqs. \ref{c42}, \ref{c43} and generate fix-point free  the 8-cell on $S^3$.
By the construction of the generators of $Q$ in {\bf Table 4.1} we also have shown the group/subgroup relation
\begin{equation}
\label{c38}
O(4,R)>G>{\rm deck}(C3)=Q.
\end{equation}

\subsection{The reduction $O(4,R)>Q$ and harmonic analysis on (C3).}

The harmonic analysis for the 3-cube $(C3)$ has as its basis 
the $Q$-periodic spherical harmonics on $S^3$ of degree $(2j)$.
Selection rules eliminate all contributions with $2j={\rm odd}$. 
The periodic  states are degenerate with respect to the subgroup 
$SU^r(2,R)$, but the degenerate states can be labelled by
a representation  index $m_2,\: -j\leq m_2\leq j$. 
The explicit construction of the $Q$-periodic polynomials 
can be found by the application of Young operators as described 
in \cite{KR08} section 4.7.

\subsection{The multiplicity $m(Q (j,j),0)$  of representations\\ in $O(4,R)>Q$.}

The reduction of representations we analyze in section 5  in the 
scheme eq. \ref{c38}. Here we  consider the reduction $O(4,R)>Q$.

To compute the multiplicity $m((j,j),0)$ for the reduction of the 
irreducible representations $D^{(j,j)}$ of $O(4,R)$ to the identity representation 
of the subgroup $Q$, denoted as $D^0(Q)$,  we need the characters $\chi^{(j,j)}$ for 
the $8$ elements of $Q$ eq. \ref{c37}. The elements $e, {\cal J}_4$ have 
the characters 
\begin{equation}
\label{c40}
\chi^{(j,j)}(e)=(2j+1)^2,\:  \chi^{(j,j)}({\cal J}_4)=(-1)^{2j}(2j+1)^2.
\end{equation}
The second equality arises because the basis of $D^{(j,j)}$ has the homogeneous degree
$2j$.
The elements $q_1,q_2,q_3$ by {\bf Table 4.3} are conjugate not in $Q$ but in $O(4,R)$ and so have the same
character w.r.t. $D^{(j,j)}$. We choose $q_1=W_2W_1W_0W_4$ as representative.
Application of \cite{KR08} eq.(60) for a product of four Weyl operators 
gives for $q_1$ in terms of the $SU^l(2,R) \times SU^r(2,R)$ action 
\begin{equation}
\label{c42}
T_{q_1}= T_{W_2W_1W_0W_4}= T_{(v_2v_1^{-1}v_0v_4^{-1},\:v_2^{-1}v_1v_0^{-1}v_4)}. 
\end{equation}
Evaluation of the two matrix products in eq. \ref{c42} with $v_0=e$ gives
\begin{equation}
\label{c43}
 (v_2v_1^{-1})(v_0v_4^{-1})=w_{l1}
=\left[ \begin{array}{ll}
0&-i\\
-i&0\\
\end{array} \right],\:
(v_2^{-1}v_1)(v_0^{-1}v_4)=w_{r1}=e.
\end{equation}
and so from eq. \ref{c43}
\begin{equation}
 \label{c44}
T_{q_1}= T_{(w_{l1},e)},\: \chi^{(j,j)}(q_1)= \chi^j(w_{l1})\chi^j(e)=\chi^j(w_{l1})(2j+1).
\end{equation}
The result eq. \ref{c43} shows that $q_1$ operates on $u\in S^3$ exclusively 
by left action corresponding to the subgroup $SU^l(2,R)< (SU^l(2,R) \times SU^r(2,R))$. The same holds true for the operators 
$T_{q_2}, T_{q_3}$ since the conjugations in eq. \ref{c35a}
applied to $T_{q_1}$ preserve the subgroup $SU^l(2,R)$.

For the characters of the element $w_l\in SU^l(2,R)$
one finds 
\begin{eqnarray}
 \label{c45}
&&\chi^{1/2}(w_{l1})= 2\cos(\phi/2)=0,
\chi^0(w_{l1})=1,
\\ \nonumber
&& \chi^{(j+2)}(w_{l1})=\chi^{j}(w_{l1})=(-1)^j
\end{eqnarray}
where the period $2$ in the last line arises from $\phi/2=\pi/2$. 

With these expressions we find the multiplicity of the $Q$-periodic states 
for a given representation $D^{(j,j)}$ as 
\begin{eqnarray}
 \label{c46}
&& m(Q (j,j),0)=\frac{1}{8}\sum_{g\in Q} \chi^{(j,j)}(g) 
\\ \nonumber 
&&= \frac{1}{8}(1+(-1)^{(2j)})(2j+1)
\left[(2j+1)+ 3(-1)^j\right].
\end{eqnarray}

The first prefactor eliminates all the states with $(2j)={\rm odd}$. 
Eq. \ref{c37} states that ${\cal J}_4 \in Q$, and therefore any $Q$-periodic state must 
be even.

The second prefactor $(2j+1)$ arises from the degeneracy with respect to the group 
$SU^r(2,R)$. For a $Q$-periodic state of polynomial degree $2j$, we can choose $(2j+1)$ 
orthogonal basis states with respect to $SU^r(2,R)$, corresponding to the second label 
$m_2$ in the spherical harmonics $D^j_{m_1,m_2}(u)$ on $S^3$.

In {\bf Table 4.2} we give the onset of the multiplicity eq. \ref{c46} for the 
lowest values of $2j$.
\vspace{0.2cm}

{\bf Table 4.1}: The three glue generators $q_i$ of ${\rm deck}(C3)$ eq. \ref{c35a} as elements
of the Coxeter group $G$ and the corresponding  pairs 
$(w_{li},w_{ri})\in SU^l(2,R)\times SU^r(2,R)$. 
\begin{equation}
\begin{array}{l|lllll}
\label{c39}
i& q_i \: x&w_{li}&w_{ri}&\epsilon&p\\
\hline
1& (x_1,-x_0,x_3,-x_2)&
\left[
\begin{array}{ll}
0&-i\\
-i&0\\
\end{array}
\right]=-{\bf k}
& e&(+-+-)&(01)(23)
\\
2& (x_2,-x_3,-x_0,x_1)&
\left[
\begin{array}{ll}
0&-1\\
1&0\\
\end{array}
\right]=-{\bf j}
& e&(+--+)&(02)(13)
\\
3& (x_3,x_2,-x_1,-x_0)&
\left[
\begin{array}{ll}
-i&0\\
0&i\\
\end{array}
\right]=-{\bf i}
& e&(++--)&(03)(12)
\\
\end{array}
\end{equation}

{\bf Table 4.2} The character $\chi^j(w_{l1})$ and the multiplicity 
$m(Q (j,j), 0),\: j \leq 8, 2j={\rm even}$, eq. \ref{c46} of $Q$-periodic spherical harmonics on $S^3$  as functions of $j$.

\begin{equation}
\label{c48}
\begin{array}{lll|lllllllll}
&\phi/2&j:&0&1&2&3&4&5&6&7&8\\
\cline{1-12}
&&&&&&&&&&\\
\chi^j(w_{l1})&\pi/2&&1&-1&1&-1&1&-1&1&-1&1\\
m(Q (j,j),0)&&&1&0&10&7&27&22&52&45&85\\
\end{array}
\end{equation}
 
\subsection{An orthogonal $Q$-periodic basis for the harmonic analysis on $C3$.}

For the quaternionic group $Q$,  the orthogonal basis for the harmonic analysis can be
given in closed analytic form. We shall use the irreducible representation matrices 
$D^j_{m_1,m_2}(u)$ given in \cite{KR08} eq. (85). Since $Q<SU^l(2,C)$ acts from the left only,
its representations  are given by $D^j_{m_1,m_2}(g),\; g \in Q$. These matrix elements are given 
in {\bf Table 4.3} for general $j$. In what follows we need them only for integer $j$.
\vspace{0.2cm}

{\bf Table 4.3}: Representation matrices for $g\in Q$ and general $j$.
\begin{equation}
\label{s1}
\begin{array}{l|l}
g & D^j_{m_1,m_2}(g), -j \leq (m_1,m_2)\leq j\\
\hline
e& \delta_{m_1,m_2}\\
q_1& \delta_{m_1,-m_2} (-1)^{j-m_1}\exp(-i\pi m_1)\\
q_2& \delta_{m_1,-m_2} (-1)^{j+m_1}\\
q_3& \delta_{m_1,m_2}  \exp(-i\pi m_1)\\
{\cal J}_4& \delta_{m_1,m_2} (-1)^{2j}
\end{array}
\end{equation}
\vspace{0.2cm}
We next construct a projection operator on the identity representation $D^0$ of $Q$:
\begin{equation}
\label{s2}
P^0:= \frac{1}{8}\sum_{g\in Q} T_g= \frac{1}{4}\left[ T_e+\sum_{i=1}^3 T_{q_i} \right],
(P^0)^2=P^0. 
\end{equation}
The representation matrix of this operator for fixed integer $j$ is obtained from {\bf Table 4.3} in the simple form
\begin{equation}
\label{s3}
(P^{0,j})_{m_1,m_1'}=\frac{1}{2}\left[1+(-1)^{m_1}\right] 
\frac{1}{2}\left[ \delta_{m_1,m_1'}+\delta_{m_1,-m_1'}(-1)^j\right]. 
\end{equation}
The projection operator eq. \ref{s3} applied to a spherical harmonic does not affect 
the column index $m_2$. 
For the projection 
it proves convenient to separate the even and odd $j$ case. 
The linearly independent, orthogonal $Q$-periodic polynomials normalized according to eq. \ref{s5}   are then found by projection 
with eq. \ref{s3}.\\

{\bf 4 Theorem}: An orthonormal basis for the harmonic analysis on the cubic spherical manifold 
$C3$ is spanned by the $Q$-periodic polynomials of degree $2j, j=0,1,2,...$ in {\bf Table 4.4}.
\vspace{0.1cm}

{\bf Table 4.4}: The $Q$-periodic orthonormal basis 
$\{ \phi^{j {\rm odd}}_{m_1,m_2}, \phi^{j {\rm even}}_{m_1,m_2}\}$ for the harmonic analysis on the 
cubic spherical manifold $C3$ in terms of spherical harmonics $D^j(u)$ on $S^3$.
\begin{eqnarray}
\label{s7}
&&j={\rm odd}, j\geq 3,\, m_1={\rm even},\, 0<m_1\leq j,\, -j\leq m_2\leq j:
\\ \nonumber
&& \phi^{j {\rm odd}}_{m_1,m_2}=\frac{\sqrt{2j+1}}{4\pi}\left[ D^j_{m_1,m_2}(u)-D^j_{-m_1,m_2}(u)\right],
\\ \nonumber
&& m(Q (j,j),0)=\frac{1}{2}(2j+1)(j-1),
\\ \nonumber 
&&j={\rm even},\, m_1=0,\,  -j\leq m_2\leq j:
\\ \nonumber
&& \phi^{j {\rm even}}_{0,m_2}=\frac{\sqrt{2j+1}}{\sqrt{8}\pi}D^j_{0,m_2}(u),
\\ \nonumber
&&j\geq 2,{\rm even},\,0<m_1\leq j,\,m_1={\rm even}:
\\ \nonumber
&& \phi^{j {\rm even}}_{m_1,m_2}=\frac{\sqrt{2j+1}}{4\pi}\left[ D^j_{m_1,m_2}(u)+D^j_{-m_1,m_2}(u)\right],
\\ \nonumber   
&&m(Q (j.j),0)= \frac{1}{2}(2j+1)(j+2). 
\end{eqnarray}
The orthogonality and normalization is  obtained by use of eq. \ref{s5}.
The multiplicities $m(Q(j,j),0)$ agree with the values derived before  in {\bf Table 4.2}.

\section{Reduction from the Coxeter group $G$.} 

\subsection{Representations of $G$.}

For the construction of irreducible representations of $G$ eq. \ref{c2} we follow Coleman \cite{COL68} in 
the notion of the little group and little co-group. For the cyclic group $C_2$ we denote the two elements 
by $\epsilon=\pm 1$ and the two irreducible representations  by 
\begin{equation}
\label{r1}
\sigma^{\pm}: \sigma^+(\epsilon)=1,\; \sigma^-: \sigma^-(\epsilon)=\epsilon.
\end{equation}
The irreducible representations of the direct product $(C_2)^4<G$ are then denoted  by 
\begin{eqnarray}
\label{r2}
&&\mu =\{\mu_1,\mu_2,\mu_3,\mu_4\}, \mu_j=\pm,\;
\\ \nonumber
&&D^{\mu}(\epsilon)= \sigma^{\mu_1}(\epsilon_0)
\sigma^{\mu_2}(\epsilon_1)
\sigma^{\mu_3}(\epsilon_2)
\sigma^{\mu_4}(\epsilon_3).
\end{eqnarray}
We choose for $D^{\mu}$  six representatives and give the little co-groups,
\begin{eqnarray}
 \label{r3}
&&D^{\mu},\: 
\\ \nonumber
&&\mu=\{++++\},\{----\},\; K=S(4),
\\ \nonumber
&&\mu=\{+++-\},\{---+\},\; K=S(3)\times S(1),
\\ \nonumber
&&\mu=\{++--\},\{--++\}, \; K=S(2)\times S(2).
\end{eqnarray}
The little co-group $K$ for the irreducible representation $D^{\mu}$  is defined 
\cite{COL68} p. 108 as the maximal subgroup $K<S(4)$ such that 
\begin{equation}
\label{r4}
h\in H: D^{\mu}(h\epsilon h^{-1})=D^{\mu}(\epsilon).
\end{equation}
The little group is the direct product  $L= (C_2)^4 \times K< G$.
Its irreducible representations are direct products of the chosen representation 
$D^{\mu}$ with the irreducible representations $D^f$ of the little co-group $H$.
Finally the irreducible representations of $G$ are the representations induced
from the little group $L<G$. This induction, denoted by the sign $\uparrow$,  is obtained  by choosing 
coset generators 
$c_j\in S(4), j=1,..,|S(4)|/|K|$
of the little co-group $K<S(4)$ and with their help constructing
\begin{equation}
\label{r6} 
D^{(\mu,f)\uparrow}_{is,jt}(\epsilon p)
=\delta(c_j^{-1}pc_i,k\in K) D^{\mu}(c_j^{-1}\epsilon c_j) D^f_{s,t}(k).
\end{equation}
The character of this representation is 
\begin{equation}
\label{r7} 
\chi^{(\mu,f)\uparrow}(\epsilon p)
=\sum_j\delta(c_j^{-1}pc_j,k\in K) D^{\mu}(c_j^{-1}\epsilon c_j) \chi^f (k).
\end{equation}
To determine these characters we must specify for each little co-group its
coset generators. The multiplicity of the identity representation of
the subgroups $H=(C_8, Q)$ is then given by 
\begin{equation}
\label{r8}
m((\mu,f)\uparrow,0)=\frac{1}{8}\sum_{g \in H} \delta(g,\epsilon p) \chi^{(\mu,f)\uparrow}(\epsilon p). 
\end{equation}
We give the corresponding data and the multiplicities in the next subsection.

\subsection{The reductions $G>C_8,\; G>Q$.}

The factorizations $g=\epsilon p$ eq. \ref{r0} of the subgroup elements are given in {\bf Tables 3.2, 4.1}. 

Representations of $G$ with $\mu=\{++++\},\{----\}, K=S(4)$: For all elements $g\in C_8, g\in Q: g=\epsilon p$ one finds in both representations $D^{\mu}(\epsilon)=1$. 
\vspace{0.3cm}
 
{\bf Table 5.1} Characters $\chi^{(\mu f)\uparrow}(g_1^t)$ for $g \in C_8$ and multiplicity 
$m(C_8)
=(C_8 \mu f),0)$ eq. \ref{r8} of the identity representation of $C_8$
in the representations 
$\mu=\{++++\}, \{----\}$
for all partitions $f$ of $S(4)$.
\begin{equation}
 \begin{array}{ll|lllll}
\chi^{(\mu f)\uparrow}(g_1^t), g\in C_8&t&[4]&[1111]&[31]&[211]&[22]\\  \hline 
&1&1&-1&-1&1&0\\
&2&1&1&-1&-1&2\\
&3&1&-1&-1&1&0\\
&4&1&1&3&3&2\\
&5&1&-1&-1&1&0\\
&6&1&1&-1&-1&2\\
&7&1&-1&-1&1&0\\
&8&1&1&3&3&2\\
\hline
&m(C_8)&1&0&0&1&1\\
 \end{array}
\end{equation}
\vspace{0.3cm}

{\bf Table 5.2} Characters $\chi^{(\mu f)\uparrow}(g)$ for $g \in Q$ and multiplicity 
$m(Q)=m((\mu f),0)$ eq.\ref{r8} of the identity representation of $Q$
in the representations 
$\mu=\{++++\}, \{----\}$. The elements $(q_1,q_2,q_3)$ are conjugate in $G$ and so their
characters coincide. 
\begin{equation}
 \begin{array}{ll|llllll}
\chi^{(\mu f)\uparrow}(g), g\in Q&g&[4]&[1111]&[31]&[211]&[22]\\  \hline 
&e&1&1&3&3&2\\
&(q_1,q_2,q_3)&1&1&-1&-1&2\\
&{\cal J}_4&1&1&3&3&2\\
&(q_1,q_2,q_3){\cal J}_4&1&1&-1&-1&2\\
\hline
&m(Q)&1&1&0&0&2\\
 \end{array}
\end{equation}
\vspace{0.3cm}

Representations with $\mu=\{+++-\},\{---+\},\: K=S(3)\times S(1)$: 
In {\bf Table 5.3} we list the  4 coset generators $c_j$  of $K=S(3)\times S(1)<S(4)$
in cycle form and give their action on $D^{\mu}(\epsilon)$.
\vspace{0.3cm}

{\bf Table 5.3} Coset generators $c_j$  and their action on $D^{\mu}(\epsilon)$ for the representations 
of $G$ with $\mu=\{+++-\},\{---+\},\: K=S(3)\times S(1)$. 

\begin{equation}
\begin{array}{l|lll}
j&c_j \in S(4)/S(3)\times S(1)& D^{\{+++-\}}(c_j^{-1}\epsilon c_j)&D^{\{---+\}}(c_j^{-1}\epsilon c_j)\\
\hline 
j&e&\epsilon_3&\epsilon_0\epsilon_1\epsilon_2\\
1&(03)&\epsilon_0&\epsilon_3\epsilon_1\epsilon_2\\
2&(13)&\epsilon_1&\epsilon_0\epsilon_3\epsilon_2\\
3&(23)&\epsilon_2&\epsilon_1\epsilon_3\epsilon_2\\
 \end{array}
\end{equation}
\vspace{0.3cm}

The representations of the little co-group $K=S(3)\times S(1)$ are 
$f= [3]\times [1],[111]\times[1],[21]\times[1]$. The characters of the elements
of $H=C_8,Q$ are given by eq. \ref{r7}. The conjugations $p\rightarrow c_j^{-1}pc_j$
cannot change the class of $p$ encoded by its cycle expression. Therefore the condition  
$c_j^{-1}gc_j=k\in S(3)\times S(1)$ eliminates all contributions except for $g=(e,{\cal J}_4)$.
The reduced set of characters and the multiplicities are given in the 
following Tables.
\vspace{0.3cm}

{\bf Table 5.4} Non-vanishing characters $\chi^{(\mu f)\uparrow}(g_1^t)$ eq. \ref{r7} for $g_1 \in C_8$ and multiplicity 
$m(C_8)=m((\mu f),0)$ eq. \ref{r8} of the identity representation of $C_8$
in the representations 
$\mu=\{+++-\}, \{---+\}$
for all partitions $f=f_1\times f_2$ of $S(3)\times S(1)$.
\begin{equation}
 \begin{array}{ll|lll}
\chi^{(\mu f)\uparrow}(g_1^t), g_1^t \in C_8&t&[3]\times [1]&[111]\times [1]&[21]\times [1]\\  \hline 
&4&-4&-4&-4\\
&8&4&4&4\\
\hline
&m(C_8)&0&0&0\\
 \end{array}
\end{equation}
\vspace{0.3cm}
 
{\bf Table 5.5} Non-vanishing characters $\chi^{(\mu f)\uparrow}(g)$ for $g \in Q$ and multiplicity 
$m(Q)=m((\mu f),0)$ eq.\ref{r8} of the identity representation of $Q$
in the representations \newline
$\mu=\{+++-\}, \{---+\}$.
\begin{equation}
 \begin{array}{ll|lll}
\chi^{(\mu f)\uparrow}(g), g\in Q&g&[3]\times [1]&[111]\times [1]&[21]\times [1]\\   \hline 
&e&4&4&4\\
&{\cal J}_4&-4&-4&-4\\
\hline
&m(Q)&0&0&0\\
 \end{array}
\end{equation}
\vspace{0.3cm}

Representations $D^{\mu}$ with $\mu=\{++--\},\{--++\},\: K=S(2)\times S(2)$: We list the 
6 coset generator of the little co-group  in the next Table.
\vspace{0.3cm}

{\bf Table 5.6} Coset generators $c_j$  and their action on $D^{\mu}(\epsilon)$ for the representations 
of $G$ with $\mu=\{++--\},\{--++\},\: K=S(2)\times S(2)$. 

\begin{equation}
\begin{array}{l|lll}
j&c_j \in S(4)/S(2)\times S(2)& D^{\{++--\}}(c_j^{-1}\epsilon c_j)&D^{\{--++\}}(c_j^{-1}\epsilon c_j)\\
\hline 
1&e&\epsilon_2\epsilon_3&\epsilon_0\epsilon_1\\
2&(12)&\epsilon_1\epsilon_3&\epsilon_0\epsilon_2\\
3&(321)&\epsilon_1\epsilon_2&\epsilon_0\epsilon_3\\
4&(120)&\epsilon_0\epsilon_3&\epsilon_1\epsilon_2\\
5&(1320)&\epsilon_0\epsilon_2&\epsilon_1\epsilon_3\\
6&(02)(13)&\epsilon_0\epsilon_1&\epsilon_2\epsilon_3\\
 \end{array}
\end{equation}
\vspace{0.3cm}

For $g_1^l \in C_8$ and the coset generators in {\bf Table 5.6} we now check the condition 
$c_j^{-1}g_1^lc_j \in S(2)\times S(2)$. The cycle structure of
$S(2)\times S(2)$ admits only the elements $e, (01),(23),(01)(23)$, and
this condition applies  at most for the elements $g_1^2,g_1^4,g_1^6,g_1^8=e$. 
For these elements one finds the non-vanishing contributions given in 
{\bf Table 5.7}.
\vspace{0.3cm}

{\bf Table 5.7} Nonvanishing character contributions for $C_8$ for $\mu=\{++--\},\{--++\}$.
The irreducible representations $f=f_1\times f_2$ of $S(2)\times S(2)$ are $D^{f_1\times f_2}((01)(23))=(-1)^{\rho_1+\rho_2},
$ with $\rho_{1,2}= 0,1$ for $f_{1,2}=[2],[11]$.
\begin{equation}
 \begin{array}{l|llllll}
g&c_j^{-1}\epsilon c_j& c_j& c_j^{-1}gc_j& D^{\{++--\}}(c_j^{-1}\epsilon c_j)
&D^{\{--++\}}(c_j^{-1}\epsilon c_j)&\chi^f(c_j^{-1}gc_j)\\ \hline
g_1^2&(-+-+)&c_3&(01)(23)& -1&-1&(-1)^{\rho_1+\rho_2}\\
     &(-+-+)&c_4&(01)(23)& -1&-1&(-1)^{\rho_1+\rho_2}\\
g_1^4&(----)&c_1...c_6&e&1&1&1\\
g_l^6&(+-+-)&c_3&(01)(23)&-1&-1&(-1)^{\rho_1+\rho_2}\\
      &(+-+-)&c_4&(01)(23)&-1&-1&(-1)^{\rho_1+\rho_2}\\
g_1^8&(++++)&c_1...c_6&e&1&1&1\\
 \end{array}
\end{equation}
Evaluation of the multiplicity eq. \ref{r8} from this Table gives for both representations 
\begin{equation}
\mu=\{++--\},\{--++\}:\: m(C_8)=m(C_8 (\mu f)\uparrow,0) = \frac{1}{8}\left[12-4(-1)^{\rho_1+\rho_2}\right]. 
\end{equation}

{\bf Table 5.8} Nonvanishing character contributions for $Q$ for $\mu=\{++--\},\{--++\}$ in
the notation of {\bf Table 5.7}.
\begin{equation}
 \begin{array}{l|llllll}
g&c_j^{-1}\epsilon c_j& c_j& c_j^{-1}pc_j& D^{\{++--\}}(c_j^{-1}\epsilon c_j)
&D^{\{--++\}}(c_j^{-1}\epsilon c_j)&\chi^f(c_j^{-1}gc_j)\\ \hline
e,{\cal J}_4&\pm(++++)&c_1...c_6&e& 1&1&1\\
q_1,q_1{\cal J}_4&\pm(+-+-)&c_1&(01)(23)&-1&-1&(-1)^{\rho_1+\rho_2}\\
   &\pm(+-+-)&c_6&(01)(23)&-1&-1&(-1)^{\rho_1+\rho_2}\\
q_2, q_2{\cal J}_4&\pm(+--+)&c_2&(01)(23)&-1&-1&(-1)^{\rho_1+\rho_2}\\
   &\pm(-++-)&c_5&(01)(23)&-1&-1&(-1)^{\rho_1+\rho_2}\\
q_3,q_3{\cal J}_4&\pm(+-+-)&c_3&(01)(23)&-1&-1&(-1)^{\rho_1+\rho_2}\\
   &\pm(+-+-)&c_4&(01)(23)&-1&-1&(-1)^{\rho_1+\rho_2}\\
\end{array}
\end{equation}
Evaluation of the multiplicity eq. \ref{r8} from this Table gives for these representations 
\begin{equation}
\mu=\{++--\},\{--++\}:\: m(Q)=m(Q (\mu f)\uparrow,0) = \frac{1}{8}\left[12-12(-1)^{\rho_1+\rho_2}\right]. 
\end{equation}

{\bf Table 5.9} Representations $D^{(\mu f) \uparrow}$ of $G$, their dimensions,  and summary on multiplicities for
$C_8$- and $Q$-periodic states. 

\begin{eqnarray}
&&
\begin{array}{l|llllll}
D^{(\mu,f)\uparrow}&\mu&f&&&&\\ \hline
&\{++++\},\{----\}&[4]&[1111]&[31]&[211]&[22]\\ \hline
{\rm dim}((\mu f) \uparrow)&&1&1&3&3&2\\
m(C_8)&& 1&0&0&1&1\\
m(Q)  && 1&1&0&0&2\\ 
\hline
\end{array}
\\ \nonumber
&&
\begin{array}{l|llll}
D^{(\mu,f)\uparrow}&\mu&f=f_1\times f_2&&\\ \hline
&\{+++-\},\{---+\}&[3]\times [1]&[111]\times [1]&[21]\times [1]\\ \hline
{\rm dim}((\mu f) \uparrow)&&4&4&8\\
m(C_8)&& 0&0&0\\
m(Q)  && 0&0&0\\
\hline
\end{array}
\\ \nonumber
&&
\begin{array}{l|lllll}
D^{(\mu,f)\uparrow}&\mu&f=f_1\times f_2&&&\\ \hline
&\{++--\},\{--++\}&[2]\times [2]&[2]\times [11]&[11]\times [2]&[11]\times [11]\\ \hline
{\rm dim}((\mu f) \uparrow)&&6&6&6&6\\
m(C_8)&& 1&2&2&1\\
m(Q)  && 0&3&3&0\\
\end{array}
\end{eqnarray}
\vspace{0.3cm}

Note that $H$-periodic states arising from  different representations of $G$ are 
orthogonal.
The defining 4-dimensional representation of the Coxeter group $G$ acting on the 3-sphere appears   in this Table  
as $(\mu f)\uparrow =((+++-) [3]\times [1])\uparrow$.

The selection rules for the cubic spherical manifolds $C2,C3$ now eliminate full representations 
of the Coxeter group $G$.
The multiplicities $m(C_8)=m(Q)=0$ for the representations with $\mu=\{+++-\},\{---+\}$ of dimensions 
$\{4,8\}$ are easily understood. 
In these representations one finds for the element ${\cal J}_4$ of both subgroups from eq. \ref{r6}
\begin{equation}
D^{(\mu f)\uparrow}_{is,jt}({\cal J}_4)=(-1)\delta_{ij}\delta_{st}, 
\end{equation}
which excludes the identity representation of $C_8$ and $Q$.

\section{Alternative spherical harmonics on $S^3$.}
The authors from the Ulm group of scientists in \cite{AU05a}, \cite{AU05b}, \cite{AU08}
have carried out a similar analysis of spherical 3-manifolds and combined it with  a powerful  astrophysical analysis.
Here we briefly compare the present and their algebraic and topological concepts.
The Laplacian $\Delta$ on $S^3$ used in \cite{AU05a} equals a Casimir operator $\Lambda^2$ of $SO(4,R)$ analyzed in 
\cite{KR05} eq.(30). The index $\beta$ is given by $\beta=2j+1$. The spherical 
harmonics on $S^3$ of \cite{AU05a} eq. (9) for fixed $j$ are related to ours algebraically by
\begin{eqnarray}
\label{A1}
&&\psi^{S^3}_{\beta lm}(u)=\delta_{\beta, 2j+1} \sum_{m_1,m_2} D^j_{m_1,m_2}(u) \langle j-m_1jm_2|lm\rangle (-1)^{j-m_1},
\\ \nonumber
&&l=0,1,2,..., 2j=\beta-1,
\\ \nonumber 
&& D^j_{m_1,m_2}(u)= \delta_{\beta, 2j+1}\sum_{l=0}^{2j}, \psi^{S^3}_{\beta lm}(u)
\langle j-m_1jm_2|lm\rangle (-1)^{j-m_1},\: m=-m_1+m_2 .
\end{eqnarray}
The coefficients are the well-known Wigner coefficients of $SU(2,C)$, \cite{ED57} pp. 31-45. 
The spherical harmonics eq. \ref{A1}  are adapted to the conjugation action
$T_{(g,g)}$: From eq. \ref{A1} with \cite{KR08} eq. (44) follows  the transformation law 
\begin{equation}
 \label{A2}
(T_{(g,g)}\psi^{S^3}_{\beta lm})(u):= \psi^{S^3}_{\beta lm}(g^{-1}ug)=\sum_{m'=-l}^l \psi^{S^3}_{\beta lm'}(u)D^l_{m',m}(g),
\end{equation}
corresponding to the law for spherical harmonics $Y_{lm}$ on $S^2$.

With respect to the topology, the authors of \cite{AU05a}, \cite{AU05b}
select homogeneous spaces which in our scheme arise from the left action 
of binary subgroups of $SU^l(2,C)$ on $S^3$. For these cases the spherical harmonics $D^j_{m_1,m_2}(u)$ eq. \ref{A1} provide the simplest
transformation law and the degeneracy (2j+1), and are used in the action of $Q$ 
on the cubic 3-manifold $C3$ in section 4.

The groups considered in \cite{KR05}, \cite{KR08} and here are derived from the homotopies of 
specific 3-manifolds and corresponding homogeneous spaces. The action of $C_5$ for the tetrahedral 
3-manifold analyzed in \cite{KR08} and of $C_8$ on the cubic 3-manifold $C2$ 
are outside   the exclusive left-or-right-action schemes. They illustrate a greater variety of 3-manifolds, homotopies, and 
subgroup actions on $S^3$. Our use of characters for multiplicities and of algebraic projection operators  avoids a numerical analysis of multiplicities and basis functions.

\section{Conclusion.}
The harmonic analysis on  spherical cubic  3-manifolds  extends  
the comparative study of Platonic 3-manifolds beyond  the dodecahedron in \cite{KR05},\cite{KR06} and the 
tetrahedron in \cite{KR08}. Two distinct homotopy groups for cubic spherical  3-manifolds $C2,C3$, implicitly defined by Everitt \cite{EV04},
are identified  as the cyclic group $H= C_8$ and the quaternion group $H=Q$ respectively  and  mapped into their isomorphic
groups of deck transformations. Using the general analysis of Weyl reflection operators 
and their representations
given in \cite{KR08}, we give on $S^3$  a $C_8$-periodic orthonormal basis 
adapted to the harmonic analysis for  the cubic spherical manifolds $C2$, and another    
$Q$-periodic orthonormal basis for  $C3$.
The distinct multiplicities, the  onset, recursion and selection rules of $H$-periodic states are computed. 
When the Coxeter group $G$ is placed in between the orthogonal group 
$O(4,R)$ and $H$, entire representations are ruled out and additional orthogonalities are
provided.

\section*{Acknowledgment.}
The author thanks T. Kramer, University of Regensburg, who did substantial
algebraic computations for the Tables.

\end{document}